\documentclass{amsart}
\usepackage{graphicx}
\usepackage{amsmath, amsthm,  amsfonts, amscd, amssymb}
\usepackage{url}
\usepackage{mathrsfs}
\usepackage[all]{xy}
\usepackage{hyperref}

\DeclareMathOperator{\hol}{hol}

\DeclareMathOperator{\ad}{ad}
\DeclareMathOperator{\ev}{ev}

\DeclareMathOperator{\Bun}{Bun}
\DeclareMathOperator{\frBun}{frBun}
\DeclareMathOperator{\id}{id}
\DeclareMathOperator{\Man}{Man}

\theoremstyle{plain}
\newtheorem{theorem}{Theorem}[section]
\newtheorem{lemma}[theorem]{Lemma}
\newtheorem{proposition}[theorem]{Proposition}

\newtheorem{corollary}[theorem]{Corollary}

\theoremstyle{definition}
\newtheorem{definition}[theorem]{Definition}

\theoremstyle{remark}
\newtheorem{example}{Example}[section]

\newtheorem{remark}{Remark}[section]

\numberwithin{figure}{section}

\newcommand{\cC}{{\mathcal C}}

\newcommand{\cF}{{\mathcal F}}

\newcommand{\fg}{{\mathfrak g}}

\newcommand{\RR}{{\mathbb R}}
\newcommand{\ZZ}{{\mathbb Z}}

\renewcommand{\a}{\alpha}
\renewcommand{\b}{\beta}

\newcommand{\<}{\langle}
\renewcommand{\>}{\rangle}
\newcommand{\ds}{\displaystyle}

\begin{document}

\title[Higher String Classes]{The caloron correspondence and higher string classes for loop groups}
  \author[M.K. Murray]{Michael K. Murray}
  \address[M.K. Murray]
  {School of Mathematical Sciences\\
  University of Adelaide\\
  Adelaide, SA 5005 \\
  Australia}
  \email{michael.murray@adelaide.edu.au}

 \author[R.F. Vozzo]{Raymond F. Vozzo}
  \address[R.F. Vozzo]
  {School of Mathematical Sciences\\
  University of Adelaide\\
  Adelaide, SA 5005 \\
  Australia}
  \email{raymond.vozzo@adelaide.edu.au}
  
  \date{November 23, 2009}

\begin{abstract}
We review  the caloron correspondence between $G$-bundles on $M \times S^1$ and 
$\Omega G$-bundles on $M$, where $\Omega G$ is the space of smooth loops in the compact
Lie group $G$. We use the caloron correspondence  to define characteristic
classes for  $\Omega G$-bundles, called string classes,   by transgression of 
characteristic classes of $G$-bundles.  These generalise the string class of Killingback
to higher dimensional cohomology. 
\end{abstract}

\thanks{The authors acknowledge the support of the Australian Research Council and useful
discussions with Alan Carey, Danny Stevenson and  Mathai Varghese.  The second author 
acknowledges the support of an Australian Postgraduate Research Award.}

\subjclass[2010]{55R10, 55R35, 57R20, 81T30}

\maketitle

\tableofcontents

\section{Introduction}

The caloron correspondence was first introduced in \cite{Garland:1988} as a bijection between
isomorphism classes of $G$-instantons on $\RR^3 \times S^1$ (calorons) and $\Omega G$-monopoles on $\RR^3$, where
$\Omega G$ is the loop group of based loops in $G$.  The movitation in that case  was the study of monopoles for 
loop groups, in particular, their twistor theory. It was subsequently \cite{Garland:1989} applied to the  case of instantons 
on the four-sphere and the four-sphere minus a two-sphere and loop group monopoles on hyperbolic three-space.
The motivation for the present work however was   \cite{Murray:2003}, which used the caloron 
correspondence to relate string structures on loop group bundles and the Pontrjagyn class of $G$-bundles.
In particular it calculated an explicit de Rham representative for Killingback's string class \cite{Killingback:1987} using
bundle gerbes.  We adopt a similar approach, without using gerbes, to define higher classes
of $\Omega G$-bundles which we call string classes and discuss their properties. 

We begin in Section \ref{sec:CW-Theory} with a brief review of Chern-Weil theory for $G$-bundles and characteristic classes.
In Section \ref{S:caloron} we explain the caloron correspondence which transforms a framed $G$-bundle over $M \times S^1$ to
an $\Omega G$-bundle on $M$ and vice-versa. We show that this is an equivalence of categories between 
the category of framed $G$-bundles over manifolds of the form $M \times S^1$ and the category of $\Omega G$-bundles. 
When the framed $G$-bundle has a (framed) connection the appropriate objects to consider on the $\Omega G$-bundle are a connection and 
a suitably defined Higgs field \cite{Murray:2003} which are introduced in Section \ref{sec:Higgs-and-cc}.

In Section \ref{S:HFs and CCs} we present our main results,   we define (higher) string classes of 
$\Omega G$-bundles,  show that they 
are characteristic classes and give an explicit formula for them.   This is a generalisation of  \cite{Murray:2003}.   A central role in this discussion 
is played by the path-fibration $PG \to G$
which is well-known to be the universal $\Omega G$-bundle.  Less well-known, but also important,
is the corresponding $G$-bundle on $G \times S^1$ introduced in \cite{CarJohMur} which plays the role 
of a universal bundle for $G$-bundles over spaces of the form $M \times S^1$.   We use the Higgs 
field to construct for any $\Omega G$-bundle a classifying map generalising results of \cite{Carey:1991}.

Throughout this paper, $G$ will be a compact,  connected Lie group and all 
cohomology groups will use real coefficients.

\section{Characteristic Classes and Chern-Weil Theory}
\label{sec:CW-Theory}

\subsection{Classifying maps and characteristic classes}

In the interests of being self-contained we shall begin by giving a short overview of the theory of classifying maps and characteristic classes before moving on to the specific case we are interested in. For details see the standard texts such as \cite{Hus}. Recall  that there is a {\em universal $G$-bundle}  $EG \to BG$   with the property that
for any  $G$-bundle over $M$ there is a so-called {\em classifying map} $f \colon M \to BG$ such that $P$ is isomorphic to the pull back of $EG \to BG$ by $f$. 
The homotopy class of the classifying map is uniquely determined by the bundle and this construction establishes a bijection between  isomorphism classes of $G$-bundles
on $M$ and homotopy classes of maps from $M$ to $BG$.  The universal bundle is characterised (up to homotopy equivalence) by the fact that it is a principal $G$-bundle and that $EG$ is a contractible space. 

A {\em characteristic class}  associates to a $G$-bundle $P\to M$ a class $c(P)$ in $H^*(M)$ which is natural with respect  to 
pulling back  in the sense that if $g \colon N\to M$ is a smooth map then 
$
c(g^*P) = g^* c(P).
$
Since all $G$-bundles are pulled back from the universal bundle $EG \to BG$ and  homotopic maps induce equal maps on cohomology,  we conclude that characteristic classes are in bijective correspondence with elements of the cohomology group $H^*(BG)$.

\subsection{The Chern-Weil homomorphism}\label{SS:Chern-Weil}
One method of constructing characteristic classes is  Chern-Weil theory. Denote by $I^k(\fg)$ the algebra of all multilinear, symmetric, $\ad$-invariant functions of 
degree $k$ on $\fg$. Elements of $I^k(\fg)$ are called invariant polynomials. If $A$ is a connection on a $G$-bundle $P \to M$ with curvature
$F$ then the $2k$-form $cw_f(A) = f(F, \dots, F)$ descends to $M$ and defines a $2k$-form which we denote by the same symbol. We have the well-known

\begin{theorem}[Chern-Weil homomorphism]\label{T:Chern-Weil}
Let $P \to M$ be a $G$-bundle with connection $A$ and curvature $F$ and let $f$ be an invariant polynomial of degree $k$ on $\fg$. Then the form $f(F, \dots, F)$ on $P$ descends to a $2k$-form on $M$ which is  closed and whose  de Rham class is independent of the choice of $A$. 
\end{theorem}

We denote the form on $M$ by $cw_f(A)$,  its de Rham class  by $cw_f(P) \in H^{2k}(M)$ and the Chern-Weil homomorphism $ f \mapsto cw_f(P)$ by 
$$
cw(P) \colon I^k(\fg) \to H^{2k}(M).
$$
Notice that it follows easily from the construction that if $\psi \colon N \to M$ and 
we endow $\psi^*P \to N$ with the pullback connection $\psi^*A $ whose curvature is  $\psi^*F$ then we  have
$$
\psi^*( cw_f(A) ) = cw_f( \psi^*A)
$$
and thus
$$
\psi^*( cw_f(P)) = cw_f(\psi^*P)
$$
so that $P \mapsto cw_f(P)$ is a characteristic class. In fact if  $G$ is compact  the Chern-Weil homomorphism 
$$
cw(EG) \colon I^k(\fg) \to H^{2k}(BG)
$$
 is an  isomorphism \cite{Dup} which extends to an algebra isomorphism  
$$
cw(EG) \colon I^*(\fg) \to H^*(BG).
$$
The proof of these results can be found in many standard places such as  \cite{Dup, Kobayashi:1969} and we will not repeat them here.
However as  a final remark we shall record here an important result about invariant polynomials which we will use in section \ref{SS:higher string classes} when we examine an analogue of the Chern-Weil homomorphism for loop groups. The derivative of the $\ad$-invariance condition on a polynomial $f \in I^k(\fg)$ gives
\begin{lemma}\label{L:ad-invariance}
Let $f \in I^k(\fg)$ and  $\a_1, \ldots, \a_k$  be $\fg$-valued forms of degree $q_1, \ldots, q_k$ respectively. Then if $A$ is a $\fg$-valued $p$-form, we have
\begin{multline*}
f([\a_1, A], \a_2, \ldots, \a_k)\\
= f(\a_1, [A, \a_2], \ldots, \a_k) + (-1)^{p q_2} f(\a_1, \a_2, [A, \a_3],\ldots, \a_k) +\ldots\\
 \ldots +  (-1)^{p(q_2 + \ldots + q_{k-1})} f(\a_1, \ldots,\a_{k-1}, [A,\a_k]).
\end{multline*}
\end{lemma}

\section{The Caloron Correspondence}
\label{S:caloron}

Before discussing the caloron correspondence  we need some definitions.

\subsection{Looping bundles}

We regard the circle $S^1$ as the quotient $\RR / 2\pi  \ZZ$ and hence denote the identity in $S^1$
as $0$.  Let $\Omega G$ be the group of all smooth maps from $S^1$ into $G$ whose value at $0$ is the identity.

 \begin{definition}
Let $P \to X$ be a $G$-bundle and $X_0 \subset X$ a submanifold. We say that $P$ is {\em framed} (over $X_0$) if we have chosen  a section  $s_0 \in \Gamma(X_0, P)$.  Denote by $P_0 \subset P$ the image of $s_0$. 
\end{definition}

The two particular examples of this general concept that we will be using are:

\begin{example} If $X$ is a pointed space, that is a space with a point $x_0 \in X$ chosen, then a framed bundle is a {\em pointed} bundle, that is a bundle $Q \to X$ with a point $q_0$ chosen in the fibre over $x_0$.
\end{example}

\begin{example} 
When  $\widetilde P \to M \times S^1$ is a $G$-bundle over $M \times S^1$ we will always
frame with respect to the submanifold $M_0 =  M \times \{ 0 \}$.
\end{example}

 If $P \to X$ is a framed bundle let $\Omega_{P_0}(P)$ be all smooth maps
 from $S^1$ into $P$ whose value at $0$ is in $P_0 = s_0(X_0)$ and similarly let $\Omega_{X_0}(X)$ be all smooth maps from $S^1$ into $X$ whose value at $0$ is in $X_0$. Note that  $\Omega_{P_0}(P) \to \Omega_{X_0}(X)$ is an $\Omega G$-bundle. We  do not discuss the  Fr\'echet principal bundle structure of 
 $\Omega_{P_0}(P) \to \Omega_{X_0}(X)$  here, details can be found in 
\cite{Hamilton:1982}, \cite{Milnor:1984} and \cite{Pressley-Segal}. Note however that we need $G$ connected so that any $G$-bundle over the circle is trivial and thus
$\Omega_{P_0}(P) \to \Omega_{X_0}(X)$ is onto. 
We call $\Omega_{P_0}(P) \to \Omega_{X_0}(X)$ a {\em loop} bundle. Of course not all $\Omega G$-bundles over $M$ are loop bundles because not all  $M$ are of the form $\Omega_{X_0}(X)$.  However we prove below the useful fact that every $\Omega G$-bundle is the pullback of a loop bundle.

\begin{remark} 
If $ Q \to X$ is a pointed bundle then instead of $\Omega_{ \{ q_0 \} }(Q) \to 
\Omega_{\{m_0 \}}(M)$ we use the notation $\Omega(Q) \to \Omega(M)$. 
\end{remark}

\subsection{The caloron correspondence}\label{SS:caloron}

To talk about the caloron correspondence it is useful to introduce some categories. 

  Let $\Bun_{\Omega G}$ be the category whose objects are $\Omega G$-bundles and morphisms are $\Omega G$-bundle maps and let $\Bun_{\Omega G}(M)$ be the groupoid 
of all $\Omega G$-bundles $P \to M$  with morphisms those bundle maps covering the identity map on $M$, that is the group of  gauge
transformations of $P$. 

 Let $\frBun_{G}$ be the category of all framed $G$-bundles $\widetilde P \to M \times S^1$ with morphisms
only those $G$-bundle maps which preserve the framing and cover a map $N \times S^1  \to M \times S^1$ of the form $f \times \id $
for some $f \colon N \to M$. Note that such a map  sends $N_0 = N \times \{ 0 \} $ to $M_0 = M \times \{ 0 \}$. 
 In both cases there are projection functors $\Pi$ to the 
category of manifolds $\Man$ defined by $\Pi(P \to M) = M$ and $\Pi(\widetilde P \to M \times S^1) = M$ and 
in the obvious way on morphisms. 

Define a  map $\eta \colon M\to \Omega_{M_0}(M\times S^1)$  by  $\eta(m)(\theta) = (m, \theta)$. Notice that $\eta(m)(0) 
= (m, 0) \in M_0$ so this is well-defined.  If $\widetilde{P} \to M \times S^1$ is a framed $G$-bundle 
then $\Omega_{\widetilde P_0}(\widetilde{P}) \to \Omega_{M_0}(M \times S^1)$ is an $\Omega G$-bundle and we can pull it back with $\eta$ to form an $\Omega G$-bundle $\cF(\widetilde P)  = \eta^*(\Omega_{\widetilde P_0}(\widetilde{P})) \to M$. It is straightforward to check that this defines a functor
$$
\cF \colon \frBun_G    \to \Bun_{\Omega G} 
$$
which commutes with the projection functor. We will show that this functor  is an equivalence of categories. This means we can find a functor 
$$
\cC \colon    \Bun_{\Omega G} \to \frBun_G
$$
and natural isomorphisms
$$
\a \colon  \cF \circ \cC \cong\id_{\Bun_{\Omega G}}  \qquad \text{and}\qquad \b \colon \cC \circ \cF \cong \id_{\frBun_G}.
$$
We shall call the functor $\cC$ the {\em caloron transform}. For simplicity, from now on we shall write $\cC^{-1}$ for $\cF$ and call it the 
{\em inverse caloron transform}. Note however that strictly speaking it is only a pseudo-inverse, that is it inverts $\cC$ only  up to natural isomorphisms.
 
To construct the  caloron transform we follow \cite{Murray:2003}. Suppose we have an $\Omega G$-bundle $P \to M$ and consider the $\Omega G$-bundle $P \times S^1 \to M \times S^1$ where the $\Omega G$ action is trivial on the $S^1$ factor. Then use the evaluation map $\ev \colon \Omega G\times S^1 \to G$ to form the associated $G$-bundle $\widetilde{P} \to M\times S^1$. That is, define $\widetilde{P}$ by
$$
\widetilde{P} = (P\times G\times S^1)/\Omega G
$$
where $\Omega G$ acts on $P\times G \times S^1$ by $(p,g,\theta)\gamma = (p\gamma, \gamma(\theta)^{-1} g, \theta)$. Then there is a right $G$ action on $\widetilde{P}$ given by $[p, g, \theta]h = [p, gh, \theta]$ (where square brackets denote equivalence classes) and a projection $\tilde{\pi} \colon \widetilde{P} \to M\times S^1$ given by $\tilde{\pi}([p, g, \theta]) = (\pi(p), \theta)$. This action is free and transitive on the fibres (which are the orbits of the $G$ action) and  $\widetilde{P} \to M \times S^1$ is a principal $G$-bundle.  Notice that over $M_0 = M \times \{0\}$ we have the well-defined framing
$$
s_0(m, 0) = [p, 1, 0]
$$
where $p$ is any point in the fibre of $P$ over $m$ and $1 \in G$ is the identity.   We define $\widetilde P$ with this framing over $M_0$  to be $\cC(P)$.

To see that the inverse caloron transform  inverts the caloron transform up to a natural isomorphism we first define 
$$
\hat\eta \colon P \to \Omega_{\widetilde P_0}(\widetilde{P}) = \Omega_{\widetilde P_0}((P\times G\times S^1)/\Omega G)
$$
by  $\hat\eta(p)(\theta) = [p, 1, \theta]$. Notice that
\begin{align*}
\hat\eta(p \gamma) (\theta) &= [p\gamma, 1, \theta]\\
                            &= [p\gamma, \gamma^{-1}(\theta) \gamma(\theta), \theta]\\
                            &= [p, \gamma(\theta), \theta] \\
                            & = (\hat\eta(p) \gamma)(\theta).
                            \end{align*}
 Thus we have the bundle map
 \begin{equation*}
\label{diagram}
\xymatrix{
 P  \ar[r]^-{\hat{\eta}} \ar[d]& \Omega_{\widetilde P_0}(\widetilde{P}) \ar[d]\\
M \ar[r]^-{\eta}& \Omega_{M_0}(M \times S^1)
}
\end{equation*} 
which defines an isomorphism $\a_P \colon \cC^{-1}(\cC(P)) = \cC^{-1}(\widetilde{P}) = \eta^*(\Omega_{\widetilde{P}_0}(\widetilde{P})) \simeq P$ and
the natural isomorphism $ \a \colon \cC^{-1} \circ \cC \cong\id_{\Bun_{\Omega G}}$. 

For the second isomorphism we start with a $G$-bundle $\widetilde{P} \to M \times S^1$ and note that the construction of $P = \cC^{-1}(\widetilde{P})$ is such that the
fibre of $P$ over $m$ is given by 
$$
P_m= \{f \colon S^1 \to \widetilde{P} \,|\, \tilde{\pi}(f(\theta)) = (m,\theta), f(0) \in \widetilde P_0\}.
$$
The $\Omega G$ action is the pointwise action of a loop.  The fibre of $\cC(P)$ over $(m, \theta) \in M \times S^1$ is given by 
$$
(P_m \times G\times \{\theta \} )/\Omega G = \{f \colon S^1 \to \widetilde{P} \mid  \tilde{\pi}(f(\theta)) = (m,\theta)\} \times G \times \{\theta\}
$$
and we have the obvious map
$$
\b_{\widetilde P} \colon [f, g, \theta ] \mapsto f(\theta)g \in \widetilde{P}_{(m, \theta)}
$$
which is a well-defined  isomorphism of $G$-bundles and  the required natural isomorphism. We need to check
that this preserves the framing. Let $\widetilde P_0$ be the framing of $\widetilde P \to M \times S^1$. The framing 
of $\cC(P)$ over a point $(m, 0)$ is given by $[f, 1, 0]$ which maps under $\beta_{\widetilde P}$ 
to $f(0) \in P_0$ so that the natural isomorphism preserves the framings.

We have now proved:

\begin{proposition}[\cite{Garland:1988,Murray:2003}]
\label{P:OmegaG correspondences}
The caloron correspondence is  an equivalence of categories between $\frBun_{G}$ and $\Bun_{\Omega G}$ commuting in both cases with the projections to $\Man$.
\end{proposition}

\begin{remark} 
Note that being an equivalence means in particular that the caloron correspondence is a bijection
between isomorphism classes of framed $G$-bundles on $M \times S^1$ and $\Omega G$-bundles on $M$. Moreover it behaves naturally with respect to maps between these bundles.
\end{remark}

We also have

\begin{corollary}
Let $P \to M$ be an $\Omega G$-bundle and let  $\eta \colon M \to \Omega_{M_0}(M \times S^1)$ be defined by $\eta(m) (\theta) = (m, \theta)$. Then there exists a framed $G$-bundle $\widetilde P \to M \times S^1$ with the property that  $P$ is isomorphic to the pull back of the loop bundle $\Omega_{\widetilde P_0}(\widetilde P) \to \Omega_{M_0}(M \times S^1)$.
\end{corollary}

\begin{example} 
Let   $Q \to X$ be a pointed $G$-bundle and $\Omega(Q)  \to \Omega(X)$ its loop bundle.   The inverse caloron transform associates to $\Omega(Q)  \to \Omega(X)$ a
framed $G$-bundle over $\Omega(X) \times S^1$.  It is straightforward to see that this is the pullback of $Q$ by the evaluation 
map $\ev \colon \Omega(X) \times S^1 \to X$ which sends $(\gamma, 0 ) \mapsto \gamma(0)$. The framing is defined by noting that
$(\ev^*Q)_{(\gamma, 0)} = Q_{\gamma(0)}$ so we can define $s_0(\gamma, 0) = q_0 \in Q_x,$ the point of $Q,$ because the loop 
$\gamma$ satisfies $\gamma(0) = x_0,$ the point of $X$. 
\end{example}

\begin{example}[The path fibration]
\label{ex:path-fibration}
Following \cite{Carey:2000} let $PG$ be the space of paths in $G$, that is smooth maps  $p \colon \RR \to G$ such that $p(0)$ is the identity and $p^{-1} \partial p$ is periodic. Then this is acted on by $\Omega G$ and
$$
\xymatrix@C=4ex{\Omega G\ar[r]	& PG\ar[d]\\
					&G}
$$
is an $\Omega G$-bundle called the \emph{path fibration}, where the projection $\pi$ sends a path $p$ to its value at $2\pi$. $PG$ is contractible and so the path fibration is a model for the universal $\Omega G$-bundle and we have $B\Omega G = G$. 

The caloron transform of the path fibration must be a $G$-bundle on $G \times S^1$. This bundle has the following simple description \cite{CarJohMur}: Start with $G \times G \times \RR$ and define an action of $\ZZ$ by $(g, h, t) n = (g, g^n h, t+2\pi n)$. Denote the quotient by  $\widetilde{PG}$. Notice that it is a principal $G$-bundle with the $G$ action given by $[g, h, t] k = [g, hk, t]$ and the projection $\widetilde {PG} \to G \times S^1$ defined by $[g, h, t ] \mapsto (g, [t])$. To see that $\widetilde{PG}$ really is $\cC(PG)$ note that we can view the latter as $(PG \times G \times \RR)/(\Omega G \times \ZZ)$ where the action is given by $(p, g, t)(\gamma, n) = (p \gamma, \gamma(t)^{-1} g, t+2\pi n)$. Then an isomorphism $\cC(PG) \xrightarrow{\sim} \widetilde{PG}$ is given by $[p, g, t] \mapsto [p(2\pi), p(t)g, t]$.  Here 
we use the fact that $p \colon \RR \to G$ satisfies $p^{-1}\partial p$ is periodic if and only if it satisfies $p(t+2\pi n) = p(2\pi)^n p(t)$. 

It is well-known that $PG \to G$ is a universal $\Omega G$-bundle and we construct a classifying map for any $\Omega G$-bundle $P \to M$ below.  Being universal in categorical language 
means that $PG \to G$ is a terminal object in the category $\Bun^h_{\Omega G}$ where morphisms are replaced by homotopy classes of morphisms.  As the caloron transformations are
an equivalence of categories it follows that $\widetilde{PG} = \cC(\Omega G) \to G \times S^1$ is a terminal object in the category $\frBun^h_G$ of framed $G$-bundles $\widetilde P \to M \times S^1$
with morphisms given by homotopy classes of maps $ f \times \id_{S^1} \colon N \times S^1 \to M \times S^1 $ where the allowable homotopies are those of the form 
$H \times \id_{S^1} \colon [0, 1] \times N \times S^1 \to M \times S^1$ for  $H \colon [0, 1] \times N \to M$ is a homotopy between maps from $N$ to $M$. 
 \end{example}

\subsection{Higgs fields and the  caloron correspondence}
\label{sec:Higgs-and-cc}

Importantly for our purposes we can extend the caloron correspondence to bundles with connection. More precisely, we have a correspondence between framed $G$-bundles on $M \times S^1$ with framed connection and $\Omega G$-bundles on $M$ with 
connection and Higgs field (Definition \ref{D:Higgs field}).

\begin{definition} Let $P \to X$ be a framed bundle with framing $s_0 \in \Gamma(X_0, P)$. A framed connection is a
connection  $A$ on $P$ such that $s_0^*(A) = 0$.
\end{definition}

\begin{lemma}
Framed connections  exist on framed bundles. 
\end{lemma}

\begin{proof} 
As $X_0 \subset X$ is a submanifold we can choose an open cover $\{ U_\a \}_{\a \in I}$ such that
$U_\a \cong U_\a \cap X_0 \times V_\a$ for some open ball $V_\a$ in $\RR^d$ where $d = \dim(X) - \dim(X_0)$. 
Moreover there is a section $s_\a \colon U_\a \to P$.
On each $U_\a \cap X_0$ we can choose $g_\a \colon U_\a \cap X_0 \to G$ such that $g_\a s_\a$ takes
values in $P_0$ and we can extend $g_\a$ to all of $U_\a$ by making it constant in the $V_\a$ directions. In other words we can just
assume that $s_\a$ restricted to $U_\a \cap X_0$ takes values in $P_0$.  We can now take
the flat connection induced by each $s_\a$ and combine these with a partition of unity. The result 
is a framed connection.
\end{proof}

The next concept we need is that of a \emph{Higgs field} for an $\Omega G$-bundle.  Let $\Omega \fg$ be the Lie algebra
of all smooth maps from $S^1$ into $\fg$, the Lie algebra of $G$, whose value at $0$ is zero.
Of course, $\Omega \fg$ is the Lie algebra of $\Omega G$. 

\begin{definition}\label{D:Higgs field}
A \emph{Higgs field} for an $\Omega G$-bundle $P \to M$ is a map $\Phi \colon P \to \Omega \fg$ satisfying the (twisted) equivariance condition
$$
\Phi (p \gamma) = \ad ( \gamma^{-1}) \Phi (p) + \gamma^{-1}\partial \gamma,
$$
for $p \in P$ and $\gamma \in \Omega G$.
\end{definition}

\begin{lemma}\cite{Murray:2003}
Higgs fields exist.
\end{lemma}

\begin{example}[The path fibration]
\label{ex:path-fibration-conn}
A connection for the path fibration is given in \cite{Carey:2002}. Let $\a$ be a smooth real-valued function on $\RR$ such that $\a(t) = 0$ for $t \leq 0$ and $\a(t) = 1$ for all $t \geq 2\pi$. Then a connection in $PG$ is given by
$$
A = \Theta - \a \, \ad(p^{-1})  \pi^* \widehat{\Theta},
$$
where $\Theta$ is the (left invariant) Maurer-Cartan form on $G$ and $\widehat{\Theta}$ is the \emph{right} invariant Maurer-Cartan form. The curvature of this connection is
$$
F =  \frac{1}{2}\left(\a^2 -\a \right) \ad(p^{-1}) [ \pi^* \widehat{\Theta}, \pi^* \widehat{\Theta}].
$$
A Higgs field for $PG$ is given by
$$
\Phi(p) = p^{-1} \partial p.
$$
We call these the  {\em standard} connection and Higgs field for the path fibration.
\end{example}

\begin{proposition}\label{P:Caloron connections 1}
Let $P \to X$ be a framed bundle with framed connection. Then the loop bundle $\Omega_{P_0}(P) \to \Omega_{X_0}(X)$ 
has a connection and Higgs field.
\end{proposition}
\begin{proof}
If $A$ is the connection on $P$ then $\Phi(q) = A(\partial q)$ defines a  Higgs field 
$$
\Phi \colon \Omega_{P_0}(P) \to \Omega\fg.
$$

  As in \cite{Murray:2003}
we can use the connection to define a connection on the loop bundle by acting pointwise. 
\end{proof}

\begin{proposition}\label{P:Caloron connections 2}
If $\widetilde P \to M \times S^1$  is a framed bundle with framed connection then the $\Omega G$-bundle $\cC^{-1}(P) \to M$ 
 has a connection and Higgs field.
\end{proposition}
\begin{proof}
It suffices to note that connections and Higgs fields pull back.
\end{proof}

Suppose instead we are given an $\Omega G$-bundle $P$ with connection $A$ and Higgs field $\Phi$. Then we can define a form on $P\times G \times S^1$ by
$$
\tilde{A} = \ad(g^{-1}) A(\theta) + \Theta + \ad(g^{-1}) \Phi \, d\theta.
$$
This form descends to a form on $\widetilde{P}$ \cite{Vozzo:PhD} and the connection (also called $\tilde{A}$) is given by this equation considered as a form on $(P\times G\times S^1)/\Omega G$. It is a connection form since if $[X, g\zeta, x_\theta]$ is a vector at $[p,g,\theta],$ (so $X \in T_pP, \zeta \in \fg$ and $x_\theta \in T_\theta S^1$) then $[X, g\zeta, x_\theta]h = [X, gh\, \ad(h^{-1})\zeta , x_\theta]$ and so
\begin{align*}
\tilde{A} ([X, g\zeta, x_\theta]h)	&= \ad(h^{-1}g^{-1}) A(X)(\theta) + \ad(h^{-1})\zeta  + \ad(h^{-1}g^{-1})x\Phi(p)\\
						&= \ad(h^{-1}) \tilde{A}([X, g\zeta ,x_\theta])
\end{align*}
and further, the vertical vector at $[p,g,\theta]$ generated by $\zeta \in \fg$ is given by
\begin{align*}
V_\zeta 	&= \frac{d}{dt}_{|_0}[p,g \exp(t\zeta),\theta]\\
		&= [0, g\zeta, 0]
\end{align*}
and so $\tilde{A}(V_\zeta) = \zeta$. (Note also that in order to show this is well defined, one needs to check that it is independent of the lift of a vector in $\widetilde{P}$. That is, if $\hat{X}$ and $\hat{X}'$ are two lifts of the vector $X\in T_{[p,g , \theta]}\widetilde{P}$ to the fibre in $P\times G\times S^1$ above $[p,g,\theta],$ then $\tilde{A}(\hat{X}) = \tilde{A}(\hat{X}')$. We leave this to the reader to check \cite{Vozzo:PhD}.) Note that the connection defined above is a framed connection. To see this, recall from section \ref{SS:caloron} that $\widetilde{P}_0$ is given by equivalence classes of the form $[p,1,0].$ A tangent vector to $[p,1,0]\in\widetilde{P}_0$ is therefore of the form $[X, 0, 0],$ for $X$ a vector field along $p$. Thus, if $(v,x) \in T_{\{m\}\times \{0\}}(M \times S^1)$ and $p$ is in the fibre above $m$, we have $s_0^*(\tilde{A})_{(m,0)}(v, x) = \tilde{A}_{[p,1,0]}([X,0,0]) = A(X)(0) = 0.$ Therefore $\tilde{A}$ is a framed connection. 

We extend our earlier notation by defining $\Bun^c_{\Omega G}$ be the category  whose objects are $\Omega G$-bundles with connection and Higgs field and morphisms are $\Omega G$-bundle maps  preserving connections and Higgs fields.  Let $\frBun^c_{G}$ be the category whose objects are  framed $G$-bundles $\widetilde{P} \to M \times S^1$ with framed connections and with morphisms only those $G$-bundle maps which preserve the framing and connection and cover a map $N \times S^1  \to M \times S^1$ of the form $f \times \id $
for some $f \colon N \to M$.  We have

 \begin{proposition}[\cite{Murray:2003}]\label{P:OmegaG connection correspondences}
The caloron correspondence 
$$
\cC \colon    \Bun_{\Omega G}^c \to \frBun_G^c
$$
and
$$
\cC^{-1} \colon \frBun_G^c   \to \Bun_{\Omega G}^c
$$
is an equivalence of categories with the same natural isomorphisms as before.
\end{proposition}
\begin{proof}
It suffices to show that the natural transformations preserve the connections.

Suppose we have the $G$-bundle $\widetilde{P} \to M \times S^1$ with connection $\tilde{A}$. Applying the caloron construction twice gives us the $G$-bundle $\cC(\cC^{-1}(\widetilde{P}))$ and according to the discussion above the connection on this bundle is given by
$$
\cC(\cC^{-1}(\tilde{A}))_{[p,g,\theta]} = \ad(g^{-1})\tilde{A}_{p(\theta)} + \Theta + \ad(g^{-1}) \tilde{A}(\partial p) d\theta.
$$
Now $\cC(\cC^{-1}(\widetilde{P}))$ is isomorphic to $\widetilde{P}$ via the map $\b_{\widetilde{P}} \colon [f, g, \theta] \mapsto f(\theta)g.$ The pushforward of this on a tangent vector $[X, g\zeta, x_\theta]$ at $[p,g,\theta]$ is given by
$$
{\b_{\widetilde{P}}}_* [X, g\zeta, x_\theta] = X(\theta)g + \iota_{p(\theta)g}(\zeta) + \partial p(\theta)xg,
$$
where $\iota_p(\zeta)$ represents the fundamental vector field at $p$ generated by the Lie algebra element $\zeta$. (Note here that $X$ is a vector tangent to the point $p \in \cC^{-1}(\widetilde{P}) = \eta^*\Omega(\widetilde{P}),$ which means it is a vector field along the loop $p$.) Using this it is easy to see that $\b_{\widetilde{P}}^* \cC(\cC^{-1}(\tilde{A})) = \tilde{A}.$

Suppose, on the other hand, we had started with an $\Omega G$-bundle $P$ with connection $A$ and constructed $\cC^{-1}(\cC(P))$. Then the connection $\cC^{-1}(\cC(A))$ is given in terms of the connection on $\cC(P)$ by acting pointwise and since the isomorphism $\a_P \colon P \xrightarrow{\sim} \cC^{-1}(\cC(P))$ is essentially given by $p \mapsto (\theta \mapsto [p,1,\theta])$ we clearly see that $\a_P^* \cC^{-1}(\cC(A)) = A.$
\end{proof}

\section{Higgs Fields and Characteristic Classes}\label{S:HFs and CCs}

\subsection{Higgs fields and the string class}\label{SS:string class}
To illustrate the caloron correspondence let us briefly outline an application, that of \emph{string structures}. This will serve not only to give an example of the correspondence above but also as motivation for the next section in which our main results will be in some sense an extension of those presented below. The material in this section is taken from \cite{Murray:2003}.

String structures were introduced by Killingback as the string theory analogue of spin structures \cite{Killingback:1987}. Suppose we have an $\Omega G$-bundle $P \to M$. Since $\Omega G$ has a central extension by the circle (see, for example, \cite{Pressley-Segal} for details) we can consider the problem of lifting the structure group of $P$ to the central extension $\widehat{\Omega G}$ of $\Omega G$. Physically, this is related to the problem of defining a Dirac-Ramond operator in string theory.\footnote{In fact, the case that Killingback considered originally was that of a free loop bundle.  Here we shall concern ourselves with the more general case of a loop group bundle which is not necessarily a loop bundle but restrict our interest to the case of based loops.} Mathematically, one has an obstruction to doing this---a certain degree three cohomology class on the base of the bundle. This class is called the \emph{string class} of the bundle and we write $s(P) \in H^3(M)$. In \cite{Murray:2003} Murray and Stevenson give a formula for a de Rham representative of this class which, adapted from the case of free loops to based loops, is given by:

\begin{theorem}[\cite{Murray:2003}]
Let $P \to M$ be a principal $\Omega G$-bundle. Let $A$ be a connection on $P$ with curvature $F$ and let $\Phi$ be a Higgs field for $P$. Then the string class of $P$ is 
represented in de Rham cohomology by the form
$$
-\frac{1}{4\pi^2} \int_{S^1} \< \nabla \Phi, F \> \, d\theta,
$$
where $\<\,\, ,\, \>$ is an invariant inner product on $\fg$ normalised so the longest root has length squared equal to $2$ and $\nabla \Phi = d \Phi + [A, \Phi] - \dfrac{\partial A}{\partial \theta}$.
\end{theorem}

In the case where $P \to M$ is given by loops in a $G$-bundle $Q \to X$ (so $P = \Omega (Q)$ and $M = \Omega (X)$) Killingback's result (also proved in \cite{Carey:1991} and \cite{Murray:2003}) is to relate this class to the first Pontrjagyn class of $Q$. In particular if $p_1(Q)$ is the first Pontrjagyn class of $Q$ then $s(\Omega (Q))$ is given by transgressing $p_1(Q)$ to $\Omega (X)$:
$$
s(\Omega (Q)) = \int_{S^1} \ev^* p_1(Q),
$$
where $\ev \colon \Omega (X) \times S^1 \to X$ is the evaluation map. In order to obtain an analogue of this result for $\Omega G$-bundles which are not loop bundles, since we do not have the option of transgressing the Pontrjagyn class, we use the caloron correspondence. Specifically, the string class of an $\Omega G$-bundle $P \to M$ is given by integrating over the circle the first Pontrjagyn class of the caloron transform of $P$ which is a $G$-bundle $\widetilde P \to M\times S^1$. To see this consider the caloron transform connection $\tilde{A}$ on $\widetilde{P}$ as in the previous section. A calculation shows that the curvature of this connection is given by
\begin{align*}
\tilde{F} 	&= d \tilde{A} + \tfrac12 [\tilde{A}, \tilde{A}]\\
		&= \ad(g^{-1})(F + \nabla \Phi \, d\theta),
\end{align*}
for $F$ the curvature of the connection $A$ on $P$ and $\nabla \Phi$ the covariant derivative of the Higgs field as above. The first Pontrjagyn class of $\widetilde{P}$ is then given by
\begin{align*}
p_1(\tilde{P}) 	&= - \frac{1}{8 \pi^2} \< \tilde{F}, \tilde{F} \>\\
			&= - \frac{1}{8 \pi^2} ( \< F, F \> + 2\< F, \nabla \Phi\> \, d\theta ).
\end{align*}
Hence integrating over $S^1$ yields

\begin{theorem}[\cite{Murray:2003}]
Let $P \to M$ be an $\Omega G$-bundle and $\widetilde{P} \to M \times S^1$ its caloron transform. Then the string class of $P$ is given by integrating over the circle the first Pontrjagyn class of $\widetilde{P}$. That is,
$$
s(P) = \int_{S^1} p_1(\widetilde{P}).
$$
\end{theorem}

The important point to note here is that the string class is canonically associated to a characteristic class for $G$-bundles, namely the first Pontrjagyn class. Furthermore, the string class is itself a characteristic class for the $\Omega G$-bundle $P$ (see \cite{Carey:1991,Vozzo:PhD}\footnote{In fact, this will follow from our work in section \ref{SS:higher string classes}.}). So using the caloron correspondence we have calculated a characteristic class of $P$. In the next section we shall extend this idea to higher degree classes and we will see that, in fact, it is possible to construct characteristic classes for loop group bundles for \emph{any} characteristic class for $G$-bundles.

\subsection{Classifying maps for $\Omega G$-bundles}\label{SS:classification}
Since we wish to calculate characteristic classes of loop group bundles, it seems natural to try to find a classifying theory for these bundles.

The case where the loop group bundle arises as a loop bundle is covered in \cite{Carey:1991} which considers a pointed bundle $Q \to X$: To write down the classifying map of the bundle $\pi \colon \Omega(Q) \to \Omega(X)$ choose a connection for $Q \to X$. Then take a loop $\gamma \in \Omega(Q)$ (so $\gamma(0) = q_0$) and project it down to $\pi\circ \gamma \in \Omega(X)$. Lift this back up to a horizontal path $\gamma_h$ in $Q$ starting at $q_0$. That is, $\gamma_h$ is horizontal,  $\gamma_h(0) = q_0$ and $\pi\circ \gamma = \pi \circ \gamma_h$. Then the \emph{holonomy}, $\hol(\gamma) \in PG$ is determined  by $\gamma = \gamma_h \hol(\gamma)$. This covers the usual holonomy  $\hol \colon\Omega(X)\to G$ and defines a bundle map:
$$
\xymatrix{P = \Omega(Q) \ar^-{\hol}[r] \ar@<2ex>[d] & PG \ar[d]\\
		 M = \Omega(X) \ar^-{\hol}[r]	& G}
$$
Thus $\hol$ is a classifying map for the bundle $\Omega(Q) \to \Omega(X)$.

We can extend this to the case of a general $\Omega G$-bundle $P \to  M$ as follows. Consider the $\Omega G$-bundle $P\to M$. Choose a Higgs field $\Phi \colon P \to \Omega\fg$ for $P$. The equation $\Phi(p) = g^{-1}\partial g$ for $g \in PG$ has a unique solution and we define the \emph{Higgs field holonomy}, $\hol_\Phi,$ by  $\hol_\Phi(p) = g$ where $g$ solves this equation.  Note that 
$$
\Phi(p h) = \ad(h^{-1})\Phi(p) + h^{-1}\partial h
$$
and 
$$
(gh)^{-1}\partial(gh) = \ad(h^{-1}) g^{-1}\partial g + h^{-1}\partial h,
$$
so that $\hol_\Phi (p\cdot h) = \hol_\Phi(p)h$ and hence $\hol_\Phi$ descends to a map (also called $\hol_\Phi$) from $M$ to $G$ and we have 

\begin{proposition}
\label{prop:universal}
If $P \to M$ is an $\Omega G$-bundle with connection $\Phi$ then  $\hol_\Phi \colon M \to G$ is a classifying map.
\end{proposition}

\begin{remark}
Recall our comment \label{ex:path-fibration} that the $G$-bundle $\widetilde{PG} \to G \times S^1$ was universal for framed $G$-bundles over manifolds
of the form  $M \times S^1$. Proposition \ref{prop:universal} also implies how to construct a classifying map $M \times S^1 \to G \times S^1$
for any $G$-bundle $\widetilde P \to M \times S^1$. That is pick a connection $\tilde A$ for $P \to M \times S^1$ and define $h \colon M \to G$ by 
sending $m \in M$ to the holonomy of $A$ around the loop $\theta \mapsto (m, \theta)$ computed relative to the framing.  The classifying
map is then $ h \times \id_{S^1}$.  Of course if $\widetilde P$ is the caloron transform of an $\Omega G$-bundle $P \to M$ with connection $A$ and 
Higgs field $\Phi$ we have that $h = \hol_\Phi$.
\end{remark}

A natural question arises at this point: If $Q \to M$ is a $G$-bundle with connection $A$ then we can define the holonomy of a loop $\gamma \in \Omega (Q)$. However, since the loop bundle $\Omega (Q) \to \Omega (M)$ is an $\Omega G$-bundle, we can also choose a Higgs field for it and define the Higgs field holonomy of a loop $\gamma$ in this bundle. Can we choose the Higgs field $\Phi$ such that $\hol_\Phi = \hol$? Define $\Phi$ in terms of $A$ as in the proof of Proposition \ref{P:Caloron connections 1} by
$$
\Phi(\gamma) = A(\partial \gamma).
$$
 Using $\gamma = \gamma_h \hol(\gamma),$ we find
$$
\partial \gamma =\partial \gamma_h \cdot \hol(\gamma) + \iota_{\gamma_{h}}(\hol(\gamma)^{-1}\partial \hol(\gamma)).
$$
Since $\gamma_h$ is horizontal (in the sense that all its tangent vectors are horizontal), applying the connection form $A$ gives
$$
A(\partial \gamma) = \hol(\gamma)^{-1} \partial \hol(\gamma).
$$
Therefore, $\hol_\Phi = \hol$.

Recall from Section \ref{SS:caloron} that the inverse caloron transform of a $G$-bundle $\widetilde{P}$ is given by the pullback $\eta^*\Omega_{P_0} (\widetilde{P}),$ for  $\eta \colon M \to \Omega_{M_0}(M \times S^1)$ defined by $\eta (m) (\theta) = (m, \theta)$. Proposition \ref{P:OmegaG correspondences} then implies that every $\Omega G$-bundle is (isomorphic to) the pullback of a loop bundle. Namely, the loop bundle $\Omega_{P_0} (\widetilde{P}) \to \Omega_{M_0} (M \times S^1)$ where $\widetilde{P} \to M \times S^1$ is the caloron transform of $P$. This suggests that there should be a relationship between $\hol_\Phi$ and $\hol$ in general. We have

\begin{lemma}\label{L:holphi=holeta}
Let $P \to M$ be an $\Omega G$-bundle with connection $A$ and Higgs field $\Phi,$ $\widetilde{P} \to M \times S^1$ its caloron transform and $\eta$ as above. Then $\ds\hol_\Phi = \hol \circ \eta$.
\end{lemma}

\begin{proof}
If $\tilde{A}$ is the connection form on $\widetilde{P}$ then $\tilde{\Phi}\colon \Omega_{P_0}(\widetilde{P}) \to \Omega\fg$ defined by
$$
\tilde{\Phi}(\gamma) = \tilde{A}(\partial \gamma)
$$
gives us that
$$
\hol_{\tilde{\Phi}} = \hol
$$
as above. Therefore we need only show that $\hol_{\Phi} = \hol_{\tilde{\Phi}} \circ \hat\eta,$ where $\hat\eta \colon P \to \Omega_{P_0}(\widetilde{P})$ is the bundle map which covers $\eta \colon M \to \Omega_{M_0} (M \times S^1)$.

Let $p\in P$. Consider the unique horizontal path $\hat\eta(p)_{h}$ such that
$$
\tilde{\pi}(\hat\eta(p)) = \tilde{\pi}(\hat\eta(p)_h)
$$
given by projecting $\hat\eta(p)$ to $\Omega_{M_0}(M\times S^1)$ and lifting horizontally back to $\Omega_{P_0}(\widetilde{P})$. The tangent vector to the loop $\hat\eta(p)$ at the point $\theta$ is given by the derivative $\partial \hat\eta(p)_\theta$ and since $\hat\eta(p)_h$ is horizontal we have that
$$
\tilde{A} (\partial \hat\eta(p)_{h,\theta}) = 0.
$$
Now, $\hat\eta(p)_\theta = [p,1,\theta]$, so we can explicitly calculate $\partial \hat\eta(p)_\theta:$
$$
\frac{\partial}{\partial \theta} \hat\eta(p)_\theta = [0,0,1].
$$
Recall that the connection $\tilde{A}$ is given in terms of the connection $A$ and Higgs field $\Phi$ for $P$ as
$$
\tilde{A} = \ad(g^{-1})A(\theta) + \Theta + \ad(g^{-1}) \Phi \, d\theta.
$$
Therefore, we have $\tilde{A}(\partial \hat\eta(p)) = \Phi(p)$. Or, in terms of the Higgs field for $\Omega_{P_0}(\widetilde{P}),$
$$
\Phi = \tilde{\Phi} \circ \hat\eta.
$$
As above, we have
$$
\tilde{\Phi}(\hat\eta(p)) = \hol(\hat\eta(p))^{-1} \partial \hol(\hat\eta(p)),
$$
and therefore $\hol_{\vphantom{{}^1}\Phi} = \hol_{\tilde{\Phi}} \circ \eta$.

\end{proof}

\subsection{Higher string classes}\label{SS:higher string classes}
In this section we shall present our main results (Theorem \ref{T:"Chern-Weil"}). As mentioned in the introduction we are interested in developing a method for geometrically constructing characteristic classes for $\Omega G$-bundles. We will accomplish this by passing to the corresponding $G$-bundle and then writing the result in terms of data on the original loop group bundle. Many of the calculations  in this section appear in more detail in the second author's PhD thesis \cite{Vozzo:PhD}.

\begin{definition}
If $P \to M$ is an $\Omega G$-bundle with connection $A$ and Higgs field $\Phi$ and $f \in I^k(\fg)$ we define the {\em string form} by 
$$
s_f(A, \Phi) = \int_{S^1} cw_f(\tilde A) \in \Omega^{2k-1}(M)
$$
where $\tilde A$ is the connection defined by the caloron transform on the $G$-bundle $\cC(P) \to M \times S^1$.
\end{definition}

While this is a definition we need a formula for the string form to be able to work with it.  
The Chern-Weil theory tells us that if we start with an invariant polynomial $f \in I^k(\fg)$ then the element in $H^{2k}(M\times S^1)$ that we end up with is $f(\tilde{F},\ldots,\tilde{F})$ where $\tilde{F}$ is the curvature of the $G$-bundle $\widetilde{P}$ on $M\times S^1$. Note that if we write out $f(\tilde{F},\ldots,\tilde{F})$ in terms of the curvature and Higgs field on the corresponding $\Omega G$-bundle $P \to M,$ we get
\begin{align*}
cw_f(\tilde A) &=f(\tilde{F},\ldots,\tilde{F}) \\ 
&= f(F+\nabla\Phi \, d\theta,\ldots,F+\nabla\Phi \, d\theta)\\
				&= f(F,\ldots, F) + k f(\nabla\Phi \, d\theta, F,\ldots,F)
\end{align*}
since $f$ is multilinear and symmetric and all terms with more than one $d\theta$ will vanish. From now on we will adopt the convention that whenever $f$ has repeated entries they will be ordered at the end and we will write them only once. That is, whatever appears as the last entry in $f$ is repeated however many times required to fill the remaining slots. (For example, $f(F) = f(F,\ldots, F)$ and $f(\nabla\Phi, F)\,d\theta = f(\nabla\Phi, F,\ldots, F)\,d\theta$.) So integrating this over the circle gives
$$
\int_{S^1} cw_f(\tilde{A}) = k \int_{S^1} f(\nabla\Phi,F)\,d\theta
$$
and we conclude
\begin{proposition}
If $P \to M$ is an $\Omega G$-bundle with connection $A$ and Higgs field $\Phi$ then the 
string form is given by 
$$
s_f(A, \Phi) = k \int_{S^1} f(\nabla\Phi,F)\,d\theta.
$$
\end{proposition}

We can now prove

\begin{proposition}\label{P:closed}
The string form is closed. 
\end{proposition}
\begin{proof}
This can be proved directly from the formula 
$$
s_f(A, \Phi) = k \int_{S^1} f(\nabla\Phi,F)\,d\theta.
$$ 
using the same methods as in Chern-Weil theory \cite{Kobayashi:1969}  but it is simpler to just note that $cw_f(\tilde A)$ is closed and that 
integration over the fibre commutes with the exterior derivative so that $s_f(A, \Phi)$ is closed.
\end{proof}

We can now consider the de Rham cohomology class of $s_f(A, \Phi)$ in $H^{2k-1}(M)$ and we have

\begin{proposition}\label{P:independence}
The class of the  string form is independent of the choice of the connection 
and Higgs field.
\end{proposition}

\begin{proof}
Again this can be proved directly but we can also note that if $(A, \Phi)$ and $(A', \Phi')$ 
are two connections and Higgs fields for $P \to M$ then we have two corresponding connections $\tilde A$
and $\tilde A'$ for the bundle $\cC(P)  \to M \times S^1$. We know from standard Chern-Weil theory
that $cw_f(\tilde A') = cw_f(\tilde A) + d \beta$ for a form $ \beta \in \Omega^{2k-1}(M \times S^1)$.
So we have
$$
s_f(A', \Phi') = s_f(A, \Phi) + d\int_{S^1} \beta.
$$
\end{proof}

An explicit formula for $\b$ is given in \cite[Proposition 3.2.4]{Vozzo:PhD}: Let $\a$ and $\varphi$ be the difference between the two connections and Higgs fields, respectively. So $\a = A' - A$ and $\varphi = \Phi' - \Phi$. Then, since the space of connections is an affine space (and the same is true for Higgs fields \cite{Murray:2003}), we can define a one parameter family of connections and Higgs fields by
$$
A_t = A + t\a, \qquad \Phi_t = \Phi + t\varphi
$$
for $t \in [0, 1]$. We let $\tilde{\a} = \ad(g^{-1})(\a + \varphi \,d\theta)$. Now consider the corresponding connection $\tilde{A}_t$ on $\cC(P)$. If $\tilde{F}_t$ is the curvature of this connection, then a calculation shows that
$$
s_f(A', \Phi') = s_f(A, \Phi) + d \left( k \int_{0}^{1} f(\tilde{\a}, \tilde{F}_t)\, dt\right).
$$

We now define

\begin{definition} 
If $P \to M$ is an $\Omega G$-bundle and $f \in I^k(\fg)$ we define the {\em string class}
of $P$,  $s_f(P) \in H^{2k}(M)$, to be the de Rham class of $s_f(A, \Phi)$ for any 
choice of connection and Higgs field.
\end{definition}

It follows immediately from the formula 
$$
s_f(A, \Phi) = k \int_{S^1} f(\nabla\Phi,F)\,d\theta
$$
that if $\psi \colon  N \to M$ then 
\begin{align*}
\psi^*(s_f(A, \Phi)) &= k \int_{S^1} \psi^*(f(\nabla\Phi,F))\,d\theta\\
               &= k \int_{S^1} f(\psi^*(\nabla\Phi),\psi^*(F))\,d\theta\\
                &= s_f(\psi^*(A), \psi^*(\Phi))
\end{align*}
and we have
\begin{proposition}
Both the string form and the string class are natural with respect to pulling back $\Omega G$-bundles with connection and Higgs field.  In particular,
the string class defines a characteristic class for $\Omega G$-bundles. 
\end{proposition}

Recall that in Example \ref{ex:path-fibration-conn} we defined a connection $A$ and Higgs field $\Phi$
on the path fibration which we called the standard connection and Higgs field. In this case we have

\begin{proposition}[\cite{Vozzo:PhD}] 
The string form of the standard connection and Higgs field of the path fibration over $G$ is 
$$
s_f(A, \Phi) =  \left( -\frac{1}{2} \right)^{k-1} \frac{k! (k - 1)!}{(2k-1)!} \, f(\Theta, [\Theta, \Theta], \ldots, [\Theta, \Theta]),
$$
where $\Theta$ is the usual left-invariant Maurer-Cartan form on $G$.  Hence the string class $s_f(PG)$, which is independent of 
the choice of connection and Higgs field is the class of $s_f(A, \Phi)$.
\end{proposition}

\begin{remark}
We remarked earlier that if $G$ is compact we have an isomorphism $I^k(\fg) \simeq H^{2k}(BG)$ given by 
$f \mapsto cw_f(EG)$. Because $EG$ is contractible we can {\em transgress } this class. 
That is, we pull it back to $EG$,   solve for $\pi^*( cw_f(EG)) = d \rho$, and  then the restriction of
$\rho$ to a fibre of $EG \to BG$ defines an element $\tau(f)$ in $H^{2k-1}(G)$ which is well known (see for example \cite{Chern:1974, Heitsch:1974})  to be the class
of the form above. That is 
$$
\tau(f) = \left( -\frac{1}{2} \right)^{k-1} \frac{k! (k - 1)!}{(2k-1)!} \, f(\Theta, [\Theta, \Theta], \ldots, [\Theta, \Theta])
$$
\end{remark}
In the following discussion we abuse notation and confuse the two maps
$$
\tau \colon I^k(\fg) \to \Omega^{2k-1}(G)
$$
and
$$
\tau \colon I^k(\fg) \to H^{2k-1}(G).
$$

Before continuing we need a simple result about integration over the fibre:

\begin{lemma}
If  $\psi \colon N \to M$ is a smooth function then pullback and integration over the 
fibre form a commuting diagram as follows:

\begin{equation*} 
\label{eq:fibre-integration} 
\xy 
(-55,7.5)*+{\Omega^q(M \times S^1) }="1"; 
(-15,7.5)*+{\Omega^q(N\times S^1)}="2"; 
(-55,-7.5)*+{\Omega^q(M)}="3"; 
(-15,-7.5)*+{\Omega^q(N)}="4"; 
{\ar^{(\psi \times \id)^*} "1";"2"};
{\ar^{\int_{S^1}} "1";"3"};
{\ar^{\int_{S^1}} "2";"4"};
{\ar^{\psi^*} "3";"4"}
\endxy
\end{equation*}

\end{lemma}
It follows that we have a commuting diagram also on cohomology.
\begin{equation*} 
\label{eq:fibre-integration} 
\xy 
(-55,7.5)*+{H^{2k}(M \times S^1) }="1"; 
(-15,7.5)*+{H^{2k}(N\times S^1)}="2"; 
(-55,-7.5)*+{H^{2k-1}(M)}="3"; 
(-15,-7.5)*+{H^{2k-1}(N)}="4"; 
{\ar^{(\psi \times \id)^*} "1";"2"};
{\ar^{\int_{S^1}} "1";"3"};
{\ar^{\int_{S^1}} "2";"4"};
{\ar^{\psi^*} "3";"4"}
\endxy
\end{equation*}
In particular if $\psi = \hol_\Phi \colon M \to G$ is the classifying map of an $\Omega G$-bundle
with connection and Higgs field $(A, \Phi)$ we have
\begin{equation*} 
\label{eq:fibre-integration} 
\xy 
(-55,7.5)*+{H^{2k}(G \times S^1) }="1"; 
(-15,7.5)*+{H^{2k}(M\times S^1)}="2"; 
(-55,-7.5)*+{H^{2k-1}(G)}="3"; 
(-15,-7.5)*+{H^{2k-1}(M)}="4"; 
{\ar^{(\hol_\Phi \times \id)^*} "1";"2"};
{\ar^{\int_{S^1}} "1";"3"};
{\ar^{\int_{S^1}} "2";"4"};
{\ar^{\hol_\Phi^*} "3";"4"}
\endxy
\end{equation*}

Composing this with the results we have already established for the path fibration we have a commuting
diagram
\begin{equation*} 
\xy 
(-95, 0)*+{I^k(\fg)} = "5";
(-55,7.5)*+{H^{2k}(G \times S^1) }="1"; 
(-15,7.5)*+{H^{2k}(M\times S^1)}="2"; 
(-55,-7.5)*+{H^{2k-1}(G)}="3"; 
(-15,-7.5)*+{H^{2k-1}(M)}="4"; 
{\ar^{(\hol_\Phi \times \id)^*} "1";"2"};
{\ar^{\int_{S^1}} "1";"3"};
{\ar^{\int_{S^1}} "2";"4"};
{\ar^{\hol_\Phi^*} "3";"4"};
{\ar^{cw({\widetilde{PG}})} "5";"1"};
{\ar_{\tau} "5";"3"}
\endxy
\end{equation*}

This gives us

\begin{theorem}\label{T:"Chern-Weil"}
If $P \to M$ is an $\Omega G$-bundle and 
$$
s(P) \colon I^k(\fg) \to H^{2k-1}(M)
$$
is the map which associates to any invariant polynomial $f$ the string class of $P$, that
is $s(P)(f) = s_f(P),$ then the following diagram commutes
\begin{equation*} 
\xy 
(-55, 7.5)*+{I^k(\fg)} = "1";
(-15,7.5)*+{H^{2k}(M\times S^1)}="2"; 
(-55,-7.5)*+{H^{2k-1}(G)}="3"; 
(-15,-7.5)*+{H^{2k-1}(M)}="4"; 
{\ar^{\tau} "1";"3"};
{\ar^{\int_{S^1}} "2";"4"};
{\ar^{\hol_\Phi^*} "3";"4"};
{\ar^{cw({\widetilde P})} "1";"2"};
{\ar^{s(P)} "1";"4"}
\endxy
\end{equation*}
\end{theorem}

Notice that although the string form is natural we would not expect the diagram in Theorem \ref{T:"Chern-Weil"} to commute
at the level of forms unless the connection and Higgs field on $P$ are the pullback of the connection and 
Higgs field on the path fibration. While it is straightforward to see that this is true for the Higgs field
it is not true for the connection.  We can however calculate what happens directly as follows.

If we start with the $G$-bundle $\widetilde{P} \to M\times S^1$ we can pull back by the evaluation map $\ev \colon [0,1] \times \Omega_{M_0} (M\times S^1) \to  M\times S^1$ to get a trivial bundle $\ev^*\widetilde{P}$ over $[0,1] \times \Omega_{M_0} (M\times S^1)$. A section is given by
$$
h \colon [0,1]\times \Omega_{M_)}(M\times S^1) \to \ev^*\widetilde{P}; \quad (t,\gamma) \mapsto  \hat{\gamma}(t),
$$
where $\hat{\gamma}$ is the horizontal lift of $\gamma$. If $\tilde{A}$ is the connection in $\widetilde{P}$ we can pull it back to $\ev^*\widetilde{P}$ and then back to $[0,1] \times \Omega_{M_0} (M\times S^1)$ to obtain
$$
\tilde{A}' := h^*\ev^*\tilde{A}.
$$
We can calculate the curvature $\tilde{F}$ of $\tilde{A}$ and pull it back by $\ev$ to $[0,1] \times \Omega_{M_0} (M\times S^1)$ and because this is a product manifold we can decompose it into parts with a $dt$ and parts without a $dt$. Under this decomposition, we have
$$
\ev^*\tilde{F} = -\frac{\partial}{\partial t} \tilde{A}'\wedge dt + \tilde{F}',
$$
where we call the component without a $dt$ $\tilde{F}'$ since if we view the form $\tilde{A}'$ for fixed $t_0$ as a connection form on $\Omega_{M_0} (M\times S^1)$ then its curvature is $\tilde{F}'$ evaluated at $t_0$. 

Now, we want to calculate $\int_{S^1} f(\tilde{F})$. The following result is straightforward and allows us to write this integral in terms of the pull back by the evaluation map

\begin{lemma}[\cite{Vozzo:PhD}]\label{L:eta* int ev* = int}
Let $\eta \colon M\to \Omega_{M_0}(M \times S^1)$ be as in section \ref{SS:caloron}. For differential $q$-forms on $M\times S^1$ we have
$$\displaystyle \eta^* \int_{S^1} \ev^* = \int_{S^1},$$
or equivalently the following diagram commutes
\begin{equation*} 
\xy 
(-55, 7.5)*+{\Omega^q(M \times S^1)} = "1";
(-15,7.5)*+{\Omega^q(\Omega_{M_0}(M \times S^1) \times S^1)}="2"; 
(-55,-7.5)*+{\Omega^q(M)}="3"; 
(-15,-7.5)*+{\Omega^q(\Omega_{M_0}(M \times S^1))}="4"; 
{\ar^{\int_{S^1}} "1";"3"};
{\ar^{\int_{S^1}} "2";"4"};
{\ar^-{\eta^*} "4";"3"};
{\ar^-{\ev^*} "1";"2"}
\endxy
\end{equation*}

\end{lemma}

Therefore for a general $\Omega G$-bundle $P\to M,$ we have
\begin{align*}
\int_{S^1} f(\tilde{F}) &= \eta^*\int_{S^1}\ev^*f(\tilde{F})\\
				&= \eta^*\int_{S^1}f(\ev^*\tilde{F}).				
\end{align*}
So we wish to calculate explicitly $\int_{S^1}f(\ev^*\tilde{F})$. To avoid many factors of $2 \pi$ we will, for this proof, regard the circle as the interval $[0,1]$ with endpoints identified. Then we can write
$$
\int_{S^1} f(\ev^*\tilde{F}) = \int_{[0,1]} f(\ev^*\tilde{F})
$$
and so we have
\begin{align*}
k\int_{S^1} f(\nabla\Phi, F)d\theta & = \eta^*\int_{S^1} f(\ev^*\tilde{F})\\
						&= \eta^*\int_{[0,1]} f(- \frac{\partial}{\partial t} \tilde{A}'\wedge dt + \tilde{F}')\\
						&= \eta^*\int_{[0,1]} f(\tilde{F}') -k \eta^*\int_{[0,1]} f(- \frac{\partial}{\partial t} \tilde{A}', \tilde{F}') dt\\
						&= -k \eta^*\int_{[0,1]} f(- \frac{\partial}{\partial t} \tilde{A}', \tilde{F}') dt.
\end{align*}
Using the formula $\tilde{F}' = d\tilde{A}' + \frac{1}{2} [\tilde{A}', \tilde{A}'],$ we can write this as:
\begin{multline*}
-k \eta^*\left\{ \int_{[0,1]} f(\partial \tilde{A}',d\tilde{A}')dt\right. \\
\left.+ (k-1)\frac{1}{2} \int_{[0,1]} f(\partial\tilde{A}', d\tilde{A}', \ldots, d\tilde{A}',[\tilde{A}',\tilde{A}'])dt +\ldots \right.\\
...+\binom{k-1}{k-2} \left(\frac{1}{2}\right)^{k-2}\int_{[0,1]} f(\partial\tilde{A}', d\tilde{A}', [\tilde{A}',\tilde{A}'])dt\\
\left.+ \left(\frac{1}{2}\right)^{k-1} \int_{[0,1]} f(\partial\tilde{A}', [\tilde{A}',\tilde{A}'])dt \right\}
\end{multline*}
where we have written $\partial \tilde{A}'$ for $\partial\tilde{A}' / \partial t$. Thus we need to work with the general term
$$
\binom{k-1}{i} \left(\frac{1}{2}\right)^{i}\int_{[0,1]} f(\partial \tilde{A}', \underbrace{d\tilde{A}', \ldots, d\tilde{A}'}_{k-i-1},  \underbrace{[\tilde{A}',\tilde{A}'], \ldots,  [\tilde{A}',\tilde{A}']}_{i})dt.
$$
To deal with these terms we shall use integration by parts  and the $\ad$-invariance of $f$ (Lemma \ref{L:ad-invariance}).
We are now in a position to prove

\begin{proposition}\label{P:s(P)=hol^*f}
$$
s_f(A, \Phi)= \hol_\Phi ^*\tau(f) + d\chi
$$
for some $(2k-2)$ form $\chi$.
\end{proposition}

\begin{proof}

To calculate the general term given above, we integrate by parts in the $\Omega_{M_0}(M\times S^1)$ and $t$ directions giving
\begin{multline*}
\int_{[0,1]}f_i dt = \int_{[0,1]} f(d\partial \tilde{A}', \tilde{A}', d\tilde{A}', \ldots, d\tilde{A}', [\tilde{A}', \tilde{A}'])dt\\
 + i \int_{[0,1]} f(\partial \tilde{A}', \tilde{A}', d\tilde{A}', \ldots, d\tilde{A}', d[\tilde{A}', \tilde{A}'], [\tilde{A}', \tilde{A}']) dt\\
 - d\int_{[0,1]} f(\partial \tilde{A}', \tilde{A}', d\tilde{A}', \ldots, d\tilde{A}', [\tilde{A}', \tilde{A}'])dt
\end{multline*}
and
\begin{multline*}
\int_{[0,1]}f_i dt = f(\tilde{A}'_1, d\tilde{A}'_1, \ldots, d\tilde{A}'_1, [\tilde{A}'_1, \tilde{A}'_1]) - f(\tilde{A}'_0, d\tilde{A}'_0, \ldots, d\tilde{A}'_0, [\tilde{A}'_0, \tilde{A}'_0])\\
 - (k-1-i)\int_{[0,1]} f(\tilde{A}', \partial d\tilde{A}', d\tilde{A}', \ldots, d\tilde{A}', [\tilde{A}', \tilde{A}'])dt\\
 - i \int_{[0,1]} f(\tilde{A}', d\tilde{A}', \ldots, d\tilde{A}', \partial [\tilde{A}', \tilde{A}'], [\tilde{A}', \tilde{A}']) dt
\end{multline*}
where we have written $f_i$ for the integrand of the general term given earlier. Combining these gives
\begin{multline*}
(k-i)\int_{[0,1]} f_i dt = f_{i,1} - f_{i,0} - i\int_{[0,1]} f(\tilde{A}', d\tilde{A}', \ldots, d\tilde{A}', \partial[\tilde{A}', \tilde{A}'], [\tilde{A}', \tilde{A}']) dt\\
+ i(k-1-i) \int_{[0,1]} f(\partial \tilde{A}', \tilde{A}', d\tilde{A}', \ldots, d\tilde{A}', d[\tilde{A}', \tilde{A}'], [\tilde{A}', \tilde{A}'])dt\\
- (k-1-i) d\int_{[0,1]} f(\partial \tilde{A}', \tilde{A}', d\tilde{A}', \ldots, d\tilde{A}', [\tilde{A}', \tilde{A}'])dt
\end{multline*}
where we have written $f_{i,1}$ and $f_{i,0}$ for $f_i$ evaluated at $t=1$ and $0$ respectively. Using $\ad$-invariance, the term on the middle line simplifies to
\begin{multline*}
\int_{[0,1]} f(\partial \tilde{A}', \tilde{A}', d\tilde{A}', \ldots, d\tilde{A}', d[\tilde{A}', \tilde{A}'], [\tilde{A}', \tilde{A}'])dt\\
\shoveleft{ = \int_{[0,1]} f(d\tilde{A}', \partial[\tilde{A}',\tilde{A}'], \tilde{A}', d\tilde{A}', \ldots, d\tilde{A}', [\tilde{A}', \tilde{A}'])dt}\\
- 2 \int_{[0,1]} f(\partial \tilde{A}', d\tilde{A}', \ldots, d\tilde{A}', [\tilde{A}', \tilde{A}'])dt\\
- (k-2-i)  \int_{[0,1]} f(\partial \tilde{A}', \tilde{A}',d\tilde{A}',  \ldots, d\tilde{A}', d[\tilde{A}', \tilde{A}'], [\tilde{A}', \tilde{A}'])dt
\end{multline*}
and so
\begin{multline*}
(k-1-i)  \int_{[0,1]} f(\partial \tilde{A}', \tilde{A}',d\tilde{A}',  \ldots, d\tilde{A}', d[\tilde{A}', \tilde{A}'], [\tilde{A}', \tilde{A}'])dt\\
= \int_{[0,1]} f(\tilde{A}', d\tilde{A}', \ldots, d\tilde{A}', \partial[\tilde{A}',\tilde{A}'],[\tilde{A}', \tilde{A}'])dt\\
- 2 \int_{[0,1]} f(\partial \tilde{A}', d\tilde{A}', \ldots, d\tilde{A}', [\tilde{A}', \tilde{A}'])dt.
\end{multline*}
Inserting this into the formula for $\int f_i dt$ gives
\begin{multline*}
(k-i)\int_{[0,1]} f_i dt = f_{i,1} - f_{i,0} -2i\int_{[0,1]} f_i dt\\
- (k-1-i) d\int_{[0,1]} f(\partial \tilde{A}', \tilde{A}', d\tilde{A}', \ldots, d\tilde{A}', [\tilde{A}', \tilde{A}'])dt
\end{multline*}
and hence
\begin{multline*}
(k+i)\int_{[0,1]} f_i dt\\
 = f_{i,1} - f_{i,0} - (k-1-i) d\int_{[0,1]} f(\partial \tilde{A}', \tilde{A}', d\tilde{A}', \ldots, d\tilde{A}', [\tilde{A}', \tilde{A}'])dt.
\end{multline*}

So we have the following expression for $s_f(A, \Phi):$
\begin{multline*}
k\int_{S^1} f(\nabla\Phi, F)d\theta\\
= -k \eta^*\left\{ \sum_{i=0}^{k-1} \binom{k-1}{i} \left(\frac{1}{2}\right)^{i} \frac{1}{k+i} \left( f_{i, 1} - f_{i, 0} -  (k-i-1) d c_i \vphantom{\tilde{f}} \right) \right\}
\end{multline*}
where $c_i$ is the last integral in the equation above (with $i$ $[\tilde{A}', \tilde{A}']$'s).

Now since $\tilde{A}'_0 = 0$ and $h(0, \gamma) = h(1, \gamma)\hol(\gamma)$ (where $h$ is the section from earlier), we have that 
$$
\tilde{A}'_0 = \ad(\hol^{-1}) \tilde{A}'_1 + \hol^{-1} d\hol
$$
and so
$$
\tilde{A}'_1 = - d\hol \hol^{-1}.
$$
Therefore we have that $f_{i,0} = 0$ and we can calculate $f_{i,1}$ in terms of $f_{0,1}$ as follows:
\begin{align*}
f_{0,1} 	&= f(\tilde{A}'_1, d\tilde{A}'_1)\\
		&= f(- d\hol \hol^{-1}, d(- d\hol \hol^{-1}))\\
		&= (-1)^{k} \left( \frac{1}{2}\right)^{k-1} \hol^* f(\Theta, [\Theta, \Theta])
\end{align*}
and in general,
\begin{align*}
f_{i,1}	&= f(\tilde{A}'_1, d\tilde{A}'_1, \ldots, d\tilde{A}'_1, [\tilde{A}'_1,\tilde{A}'_1])\\
		&= (-1)^{k-i} \left(\frac{1}{2}\right)^{k-1-i} \hol^* f(\Theta, [\Theta, \Theta])\\
		&= (-1)^{i} 2^{i} f_{0,1}
\end{align*}
using the fact that $d(-d\hol\hol^{-1}) = -\frac{1}{2} [d\hol \hol^{-1}, d\hol \hol^{-1}]$.

Therefore we have
\begin{multline*}
k\int_{S^1} f(\nabla\Phi, F)d\theta\\
=  \left(-\frac{1}{2}\right)^{k-1} k \sum_{i=0}^{k-1} \binom{k-1}{i}  \frac{(-1)^i}{k+i} \hol_\Phi^*f(\Theta, [\Theta, \Theta])\\
+ k\sum_{i=0}^{k-i}\binom{k-1}{i} \left(\frac{1}{2}\right)^i  \frac{1}{k+i}   (k-i-1) d c_i.
\end{multline*}

It turns out \cite{Sury:2004} that the coefficient above is equal to the coefficient in the definition of the transgression map: 
$$
k \sum_{i=0}^{k-1} \binom{k-1}{i}  \frac{(-1)^i}{k+i} =  \frac{k! (k - 1)!}{(2k-1)!}.
$$
So we see that the pull back of $\tau(f)$  is cohomologous to the string form. 

\end{proof}

\section{Conclusion}
There is an immediate natural generalisation of what we have done which is to replace $M \times S^1$ by an $S^1$-bundle $Y \to M$. 
In this case the caloron correspondence has been used in  \cite{Bergman:2005} in an application to string theory and in \cite{BouMat}
in an application to $T$-duality.  String classes in this case have been constructed in \cite{Vozzo:PhD} and will appear in  \cite{MurVoz}.

The results we have presented are one way of defining characteristic classes for infinite dimensional bundles by essentially using 
the caloron correspondence to avoid the infinite dimensionality. Another approach would be to deal with the infinite dimensionality
directly by extending the notion of invariant polynomials to a genuinely infinite dimensional setting.  Paycha, Rosenberg and 
collaborators have done this by using the Wodzicki residue as a trace on the Lie algebra of the group of invertible, zeroth-order
pseudo-differential operators on a vector bundle over a compact space (see for example \cite{Paycha:2004} or \cite{Paycha:2007} and references therein). In the case of a trivial, rank $n$ real vector
bundle over the circle this group contains the loop group of $O(n)$ as the subgroup of multiplication operators. The 
Wodzicki characteristic classes defined in this way vanish on bundles whose structure group has a reduction 
 to the loop group \cite{Ros}.

Finally we note two  unanswered questions. Firstly we know from \cite{Murray:2003} that the 
three-dimensional string class is the obstruction to lifting the structure group of a  loop group bundle to 
the Kac-Moody group. The geometric significance of the higher string classes is an open question. 
Secondly if we regard $\Omega G$ as based gauge transformations of a bundle over the circle, how much of our work  can be generalised to an arbitrary group of gauge transformations?

\end{document}